
\overfullrule=0pt

\input amstex

\documentstyle{amsppt}

\define\cl{$\Cal L$} \define\clsp{$\Cal L$ }

\define\emo{$ M_0$}  \define\emosp{$ M_0$ }

\define\eno{$ N_0$}  \define\enosp{$ N_0$ }
\define\enn{$ N_n$}  \define\ennsp{$ N_n$ }
\define\emr{$ M_r$}  \define\emrsp{$ M_r$ }
\define\dpr{$ D^\prime$}  \define\dprsp{$ D^\prime$ }
\define\bethr{$ B^3$}

\define\gam{$\gamma$}  \define\gamsp{$\gamma$ }

\define\nidt{\noindent}
\define\ubr{\underbar}
\define\ssk{\smallskip}
\define\msk{\medskip}

\define\smin{$\setminus$}
 
\define\ap{$\cap$}

\define\up{$\cup$}

\define\sset{$\subseteq$} 

\define\del{$\partial$}
\define\delv{$\partial_{v}$} 
\define\delh{$\partial_{h}$}

\topmatter

\title\nofrills Persistently laminar tangles \endtitle
\author  Mark Brittenham \endauthor

\leftheadtext\nofrills{Mark Brittenham}
\rightheadtext\nofrills{Persistently laminar tangles}

\affil   University of North Texas\endaffil
\address   Department of Mathematics, Box 305118, University of North Texas, 
Denton, TX 76203 \endaddress
\email   britten\@unt.edu \endemail
\thanks   Research supported in part by NSF grants \# DMS$-$9400651 and DMS$-$9704811\endthanks

\abstract
We show how to build tangles $T$ in a 3-ball with the property that 
any knot obtained by tangle sum with $T$ has a persistent lamination in its 
exterior, and therefore has property P. The construction is based on an example
of a persistent lamination in the exterior of the twist knot $6_1$, due
to Ulrich Oertel. We also show how the construction can be generalized to 
$n$-string tangles.
\endabstract

\endtopmatter

\document

\heading{\S 0 \\ Introduction}\endheading

Essential laminations have proved very useful in understanding the topology of knots 
in the 3-sphere. Constructions of 
essential laminations in knot exteriors have allowed us to see that non-trivial surgery 
on non-torus alternating knots [DR], and on most algebraic knots [Wu] yield manifolds with 
universal cover ${\bold R}^3$, for example. This can be thought of as a (very) strong form 
of Property P for these knots. They can also provide a means of detecting the underlying 
geometric structure of the 3-manifolds obtained by surgery on a knot [Br1],[Br2],[BW].

In this paper we construct persistent laminations for knots, that is, 
essential laminations in the exterior of the knot, which remain essential after any 
non-trivial Dehn filling. Our starting point is a particularly simple 
example of such an essential lamination \cl, found by Ulrich Oertel [Oe] in the 
complement of a 
twist knot (the knot $6_1$ in Rolfsen's knot tables [Ro]), in connection with 
his work on laminations with a transverse affine structure. What we show here is that 
this lamination can be associated to a rather simple tangle $T_0$. By this we mean 
two things: (1) the lamination \clsp lives in the complement of the tangle $T_0$ in 
the 3-ball \bethr; (2) if we sum $T_0$ with any other tangle $T$ to obtain a knot $K$ 
in the 3-sphere $S^3$, then \clsp is  
{\it persistent} for $K$. We call such a tangle {\it persistently laminar}. 
Being persisently laminar immediately implies, for example that every knot $K$ obtained 
by tangle sum with $T_0$ has Property P. We also show that the construction of the 
lamination \clsp can be generalized to provide many more examples of persistently 
laminar tangles.

\heading{\S 1 \\ The lamination}\endheading

Oertel's construction of the lamination \clsp begins with the branched surface $B$,
depicted in Figure 1a, 
embedded in the complement of the $6_1$ knot $K_0$. We 
have removed the knot $K_0$ in Figure 1b, to give a better view of the branched surface. 
This branched surface can be thought of as a once-punctured torus (i.e, a disk with a 
1-handle attached), with 
its boundary glued to a curve running over the 1-handle, to create a single, embedded 
branch curve for $B$ (as well as a second tube for the knot to run through); see Figure 
1c.

\ssk

\input epsf.tex

\leavevmode

\epsfxsize=4.5in
\centerline{{\epsfbox{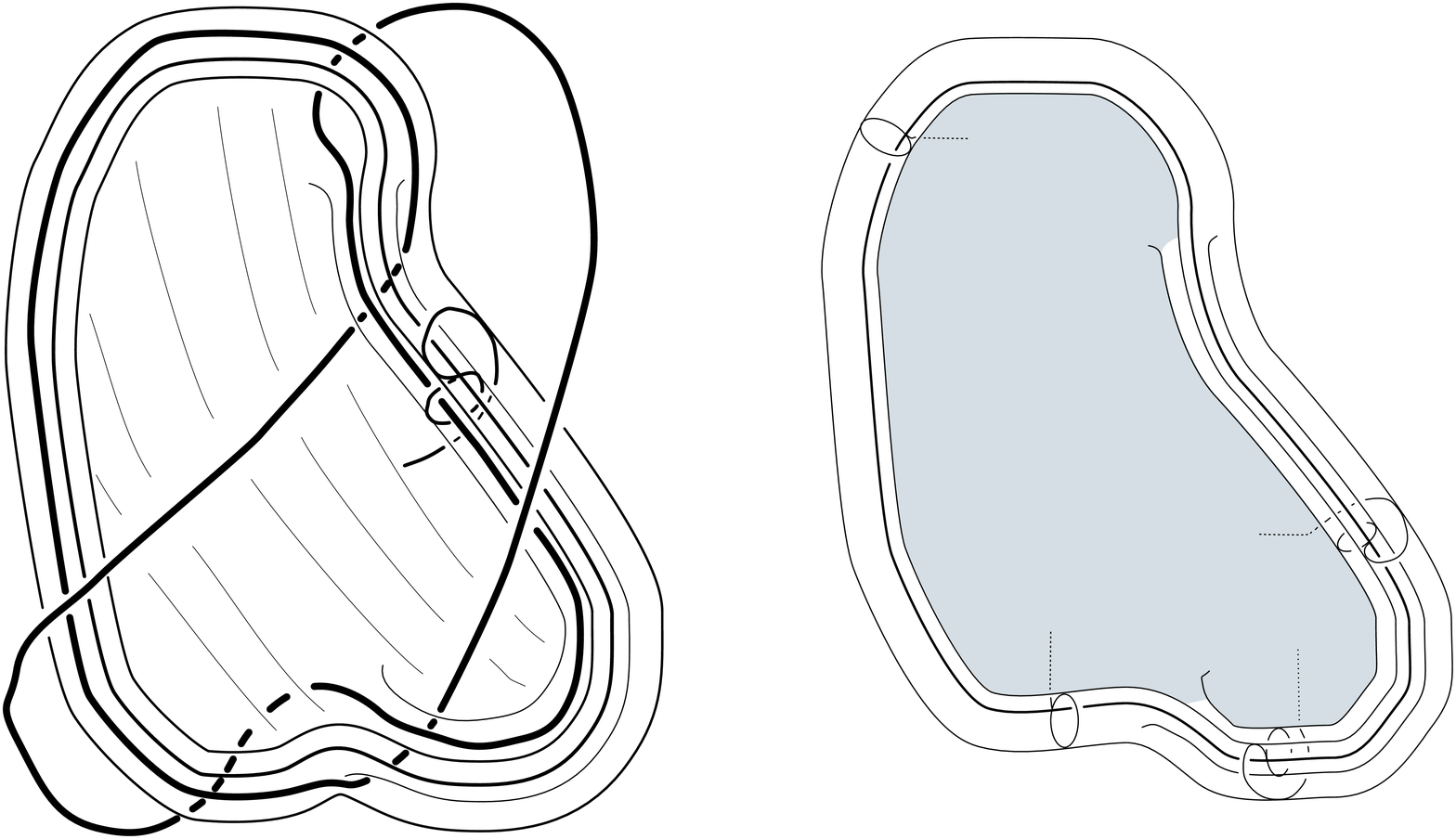}}}

\centerline{Figure 1a \hskip2in Figure 1b}

\ssk
\leavevmode

\epsfxsize=4.5in
\centerline{{\epsfbox{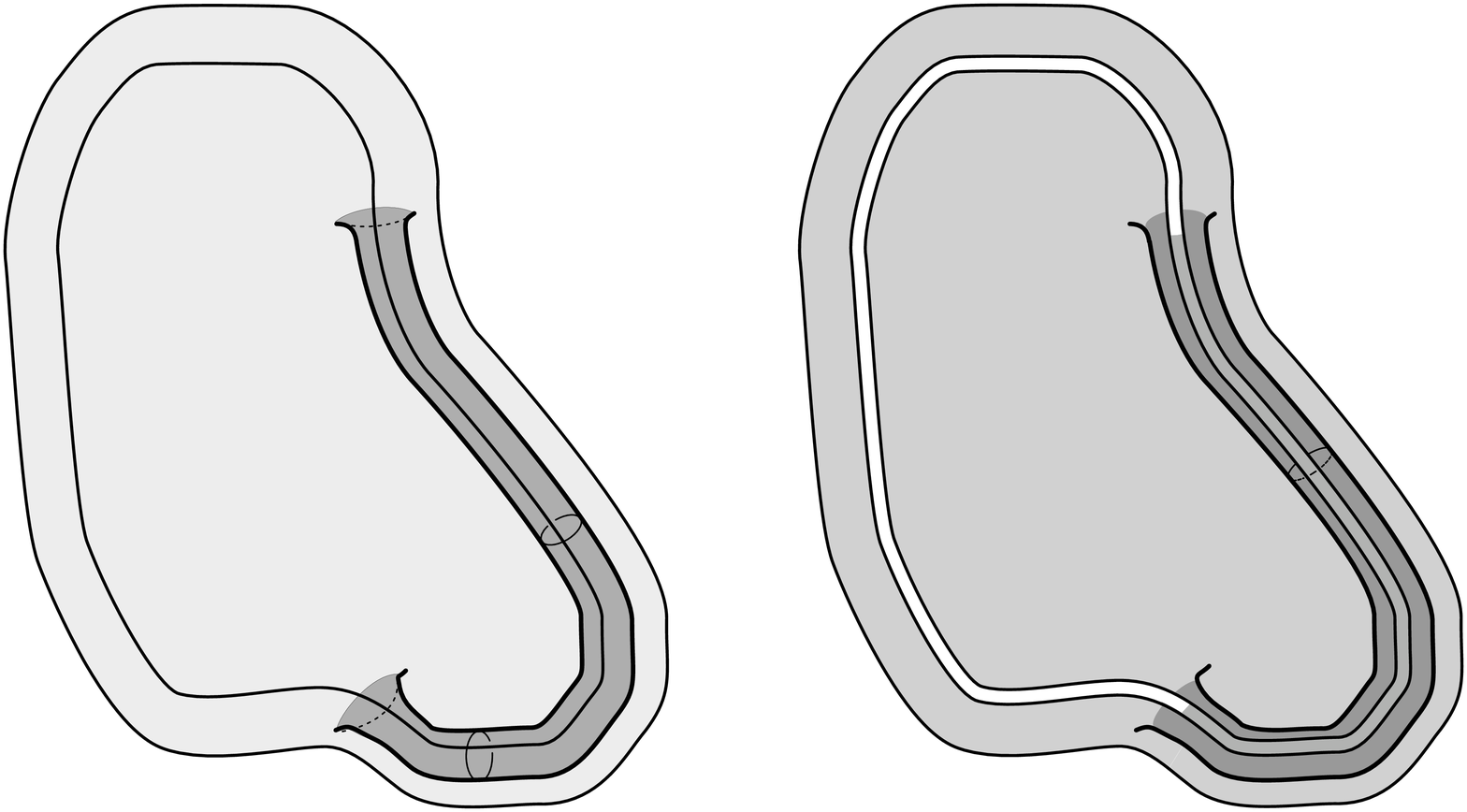}}}

\centerline{Figure 1c \hskip2in Figure 1d}

\smallskip

Figure 2 provides a magnified picture of the process of gluing at the base of the 
1-handle, to make it easier to see how the second tube is created.

Note that \enosp = $S^3$\smin int$N(B)$ is a genus-2 handlebody; this is 
most easily seen from Figures 1b and 2. By `filling in' the two tubes that the 
knot $K_0$ runs through (one of which is created when we glue the boundary of the punctured 
torus to \gam\ to create $B$), which we can think of as gluing two 2-disks 
$D_1$, $D_2$ to $B$, we can see that $N(B$\up $D_1$\up $D_2)$ is a 3-ball, 
i.e, $S^3$\smin int$N(B$\up $D_1$\up $D_2)$ is a 3-ball. \enosp is therefore 
a handlebody.

\ssk
\leavevmode

\epsfxsize=4in
\centerline{{\epsfbox{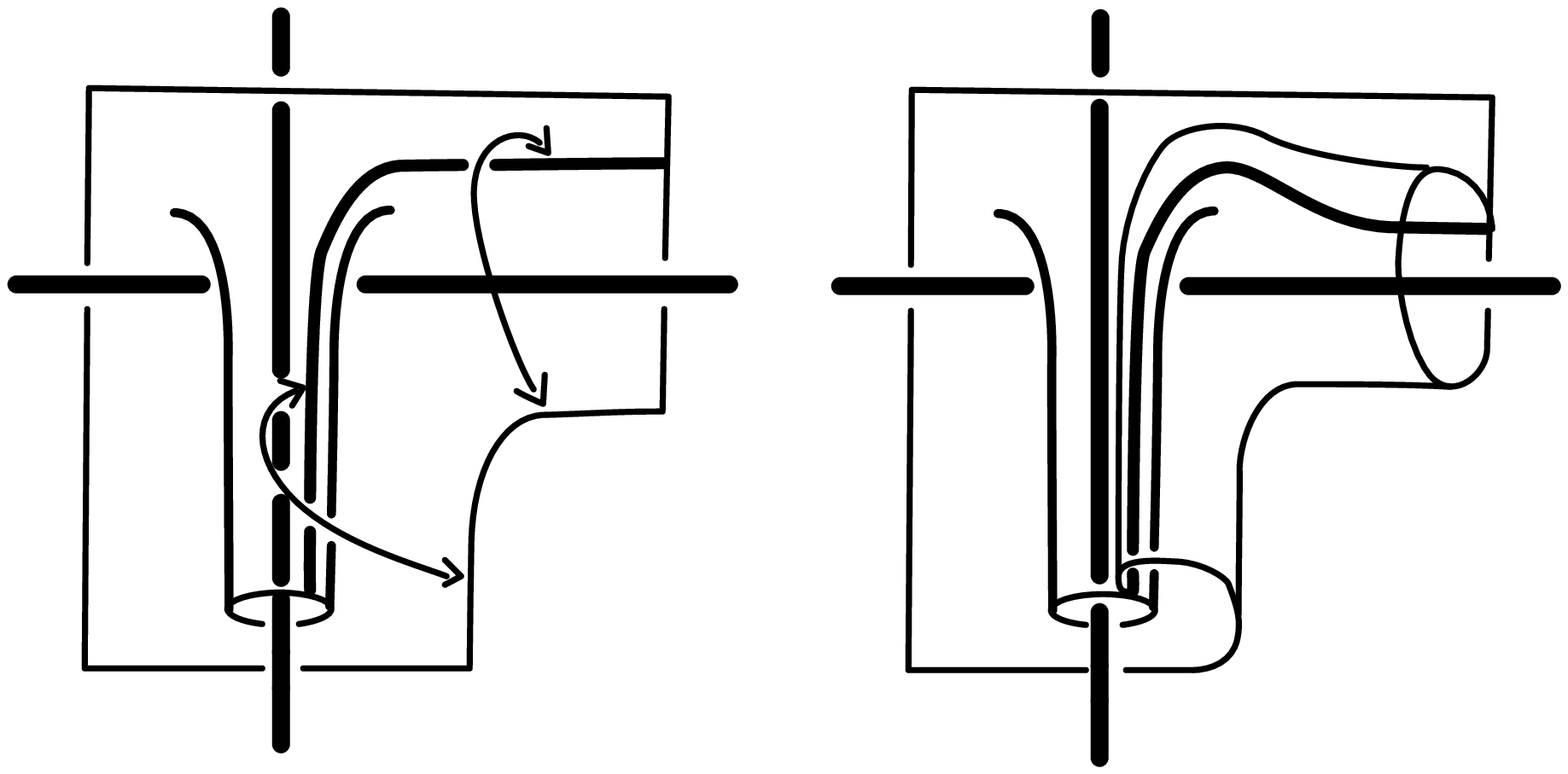}}}

\centerline{Figure 2}

\ssk

Next we show how to find a knot in the exterior of $B$ so that $B$ will be essential 
in the exterior of the knot. Our condition comes from the fact that our knot $K_0$
above meets each of the disks $D_1,D_2$ in one point.

\proclaim{Theorem} Let $K$ be a knot in \enosp = $S^3$\smin$N(B)$ meeting each of the 
disks $D_1,D_2$ in one point. Then $B$ is an essential branched surface in 
$M=S^3$\smin int($N(K)$).\endproclaim

{\bf Proof:} For $B$ to be essential we need to know 6 things:

\ssk

(1)  $B$ carries a lamination \cl\ with full support.

\ssk

\nidt This is immediate, since $B$ has no triple points; the branch curve \gamsp does 
not intersect itself. If we cut $B$ open along \gamsp (see Figure 1d), we get a surface 
with boundary, $F$. By taking a Cantors set's-worth of copies of $F$, embedded transverse 
to the fibers of $N(B)$, we can glue these surfaces together where the three copies of 
\gamsp meet (since the concatenation of two Cantor sets is order isomorphic to a Cantor 
set) to create a lamination \cl\ carried with full support by $B$. 

\ssk

(2) $B$ does not carry a 2-sphere, and $B$ has no disks of contact.

\ssk

This follows because $B$ has only one sector, i.e., $B$\smin\gamsp is connected. The 
sector is in fact a twice-punctured 2-disk $F$ (it is, after all, a once-punctured torus 
cut open along a non-separating curve); see Figure 1d. Any surface carried by $B$ 
would consist of finite number $a$ of parallel copies of $F$, glued together where 
the three boundary components of $F$ come together at \gam. When these three sheets 
come together we get a consistency condition to determine if we can glue the 
boundary components together to get a closed surface (Figure 3). In this case 
the condition is $a$ + $a$ = $a$, implying $a$ = 0. So no such surface exists. 
A disk of contact is similar; it is a 2-disk carried by $B$, whose boundary lives 
in the vertical boundary \delv $N(B)$ of $B$. This must again be built by gluing 
copies of $F$ together, except this time, after gluing, a boundary component is left 
free. This 
gives the consistency condition $a$ + $a$ +1 = $a$, implying $a$ =$-$1, which is 
absurd. So there are also no disks of contact.

In point of fact, we have 
shown that $B$ carries no closed surface, and has no compact 
surfaces of contact.

\ssk
\leavevmode

\epsfxsize=4in
\centerline{{\epsfbox{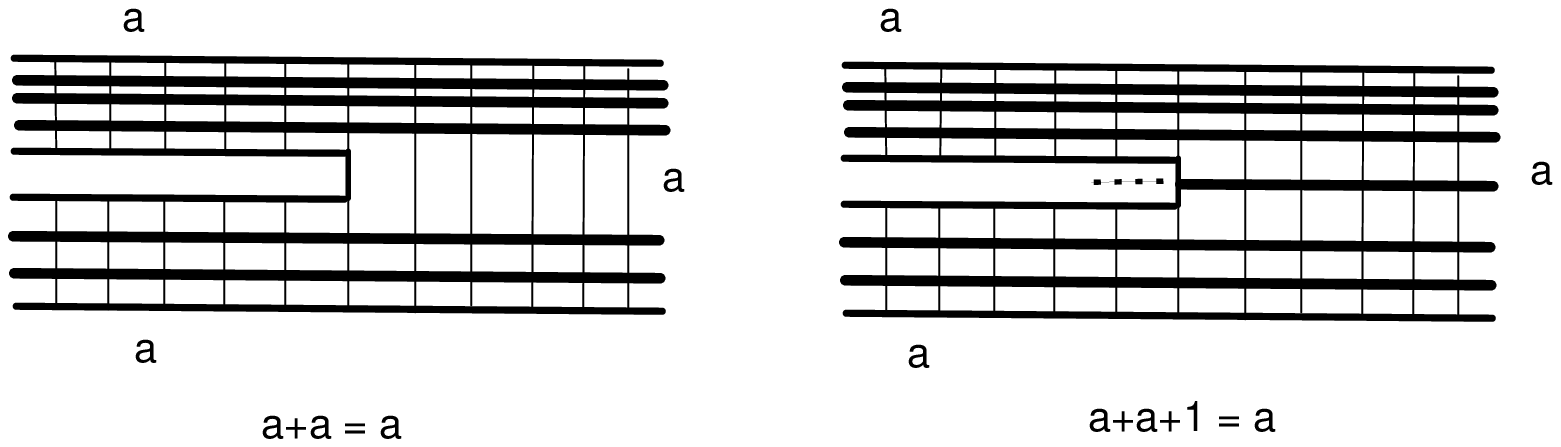}}}

\centerline{Figure 3}

\ssk

(3) $B$ does not carry a compressible torus.

\ssk

\nidt This follows from the above, since $B$ does not carry a closed surface.

\ssk

(4) \emosp = $M$\smin int($N(B)$) does not have any monogons.

\ssk

\nidt This is also immediate, because $B$ is transversely orientable; there is 
a vector field (in $M$) everywhere transverse to the tangent planes of $B$. 
The arc in the boundary of a monogon  which meets 
\delh$(N(B))$ is a transverse orientation-reversing loop.

\msk

(5) \emosp is irreducible.

\ssk

\nidt  Suppose $S$ is a reducing sphere for \emo. Since \eno\ is a handlebody, $S$
bounds a 3-ball $B^3$ in \eno. Since this 3-ball cannot live in \emo, we must have
$K$\sset$B^3$. But this implies that $K$ is null-homotopic in \eno, hence
homologically trivial. But $K$ intersects $D_1$, for example, exactly once, and
so its homology class has non-trivial intersection number with the class of $D_1$,
so is non-trivial (Figure 4). So the reducing sphere cannot exist.

\ssk

\leavevmode

\epsfxsize=2in
\centerline{{\epsfbox{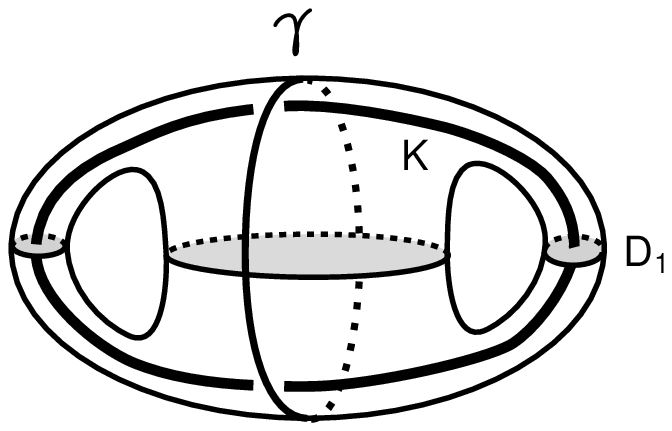}}}

\centerline{Figure 4}

\ssk

(6) The horizontal boundary \delh $N(B)$ of $B$ is incompressible in \emo.

\ssk

\nidt Again the idea is to use the fact that $K$ pierces the 
disks $D_i$ in one point each. The point is that there are very few 
compressing disks for \delh $N(B)$ 
in the handledbody \enosp to begin with, and the fact that we are 
dealing with a knot means 
that all of them must intersect $K$.

Suppose $D$ is a compressing disk for \delh $N(B)$ in \emo. Then in particular, 
it is a compressing disk for \delh $N(B)$ in \eno. The key to the argument is 
the fact that \gamsp itself bounds a compressing disk, call it \dpr. because 
\del$D$\ap\gamsp = \del$D$\ap\del\dprsp = $\emptyset$, we can, by a 
disk-swapping argument, make $D$ and \dprsp disjoint. Now \eno\smin\dprsp 
consists of two solid tori, each with a fat point (namely \dpr) removed from 
their boundaries. $D$ can therefore be thought of as a compressing disk for 
a solid torus which has a point removed from its boundary. It is therefore 
isotopic either to a meridian disk of the solid torus (and therefore, back 
in \eno, is isotopic to one of the disks $D_1$ or $D_2$), if \del$D$ is 
essential in the boundary of the solid torus, or is boundary parallel, i,e, 
is parallel to a disk in the boundary of the solid torus, with the removed 
point in its interior; see Figure 5a. Back in \enosp, this second disk is 
parallel to \dpr.

\ssk
\leavevmode

\epsfxsize=5in
\centerline{{\epsfbox{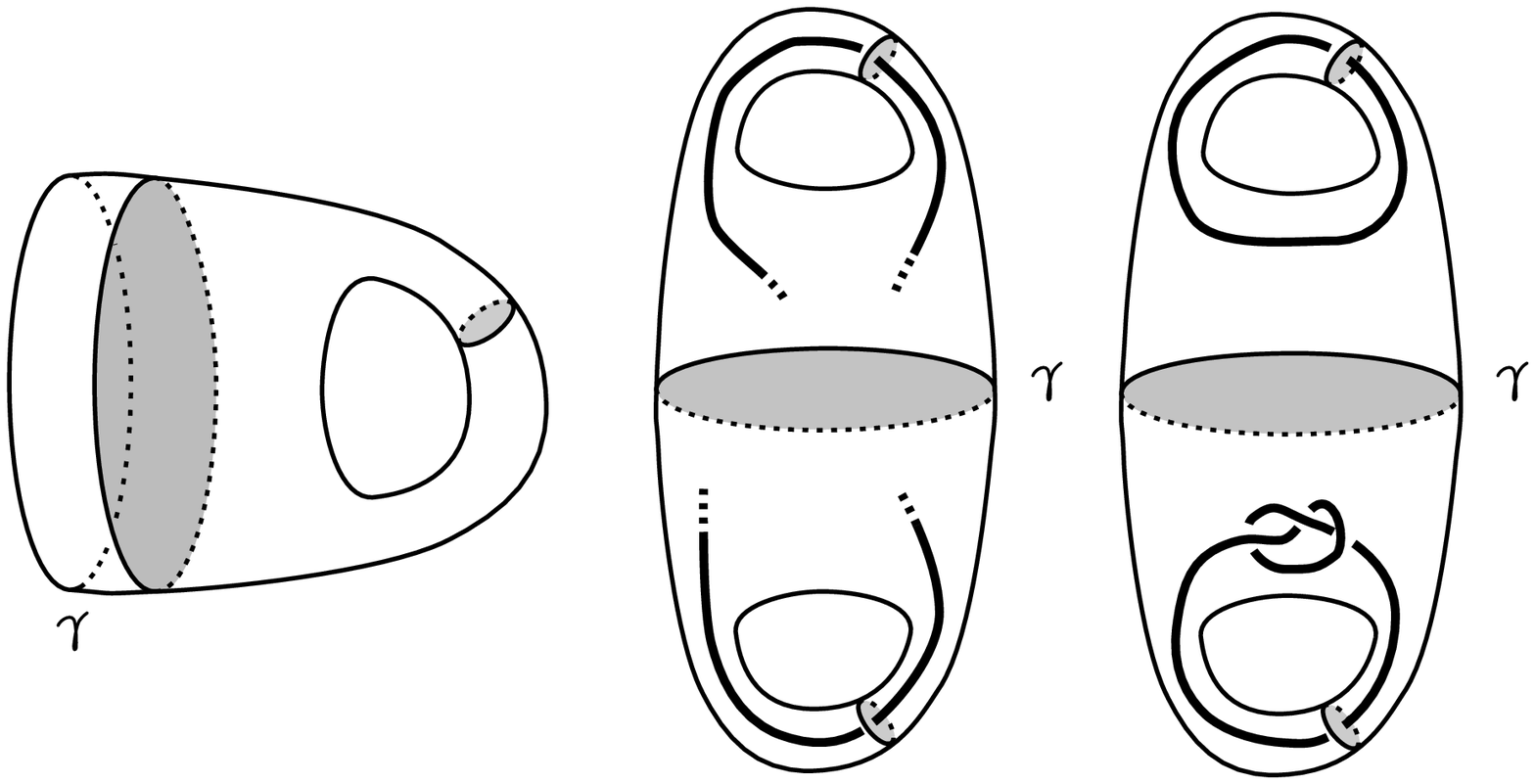}}}

\centerline{Figure 5}

\ssk

But both of these possibilities are absurd; in the first case $D$, which misses 
$K$, is isotopic, rel boundary (if we wish - the original isotopy moved \del 
$D$ around in \del\eno), to a disk which hits $K$ exactly once, contradicting 
the invariance of intersection number for homology classes (Figure 5b). In the 
second case $D$ separates \eno, yet $K$, which is connected, has non-trivial 
intersection with each piece (Figure 5c). Therefore, no compressing disk for 
\delh $N(B)$ in \emosp exists.

Consequently, all properties of essentiality are satisfied, so $B$ is an 
essential branched surface, and \clsp is an esential lamination, in $S^3$\smin 
int$(N(K))$. \,\, \qed

\ssk

\leavevmode

\epsfxsize=4.5in
\centerline{{\epsfbox{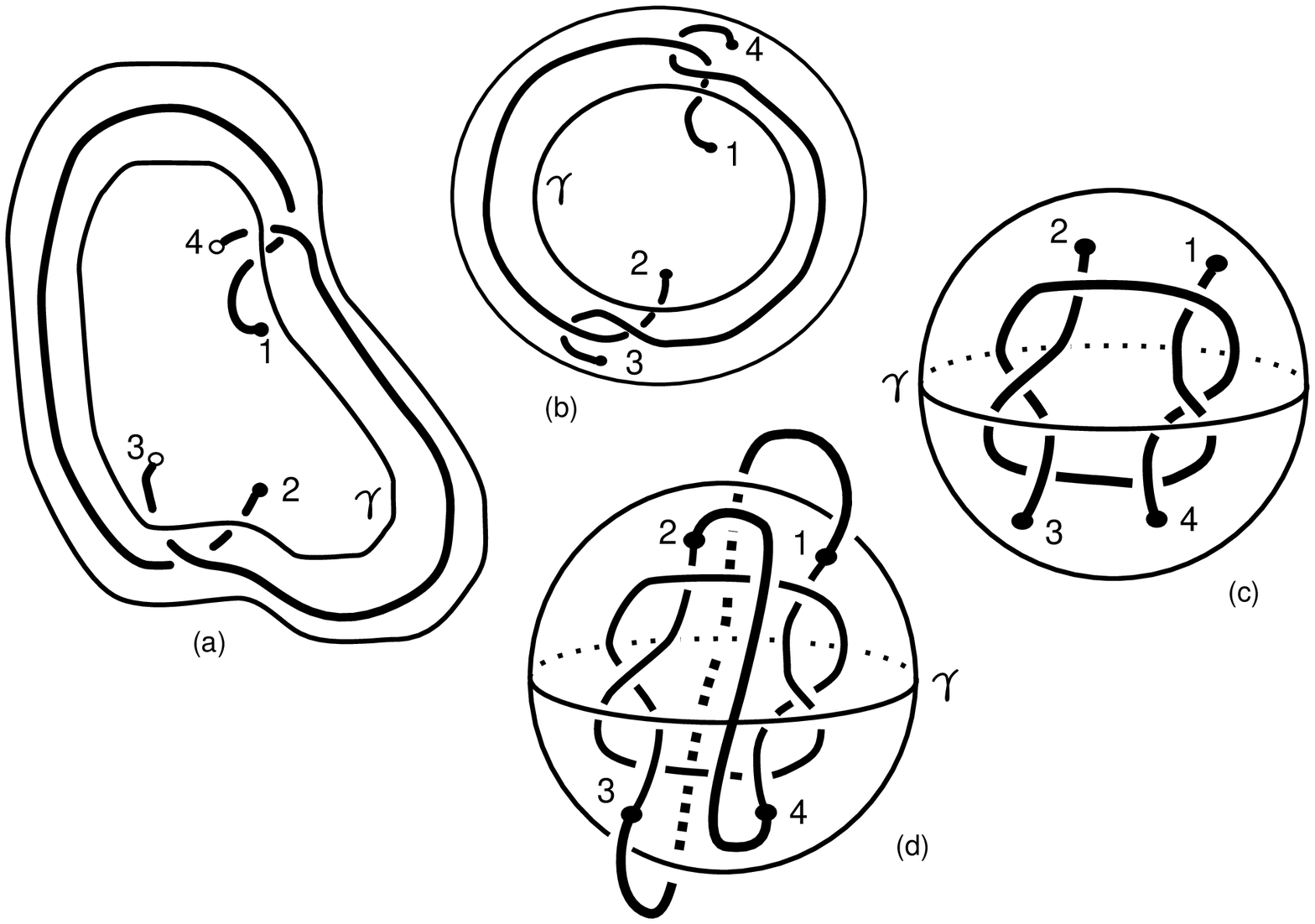}}}

\centerline{Figure 6}

\ssk

\heading{\S 2 \\ The tangle $T_0$}\endheading

$N(B$\up$D_1$\up$D_2$) is a 3-ball $B_0$, which the 
our original twist knot $K_0$, and each of the knots $K$, intersects in a 
pair of arcs, namely the two cores 
of the tubes that the disks cap off. In other words, each $K$ meets this 
3-ball in the same tangle $T_0$. But we have not yet identified this tangle. To see 
what it is, imagine stretching the two arcs in the 3-ball to fill the two 
tubes of $S^3$\smin int$(N(B))$; we arrive at a picture as in Figure 6a (we 
have also included the branch curve \gamsp for reference). The key point 
here is that the arc which, in this picture, crosses \ubr{over}  the other, 
comes out of the 3-ball in \ubr{back}. Pulling it around front adds extra 
half-twists to the tangle, so that we end up with the tangle in Figure 6c. 
It is the sum of two rational tangles, the 1/3 tangle and the $-$1/3 tangle.

Our original knot $K_0$ is obtained from this tangle by tangle sum with a rational 
tangle; see Figure 6d.

\heading{\S 3 \\ Persistence}\endheading

The analysis of Section 1 above shows that the lamination \clsp is essential in 
the complement of any knot $K$ obtained by tangle sum of $T_0$ with any other 
tangle $T$. We now show:

\proclaim{Proposition} \clsp is persistent for the knot $K$; it remains essential in 
any manifold 
obtained by non-trivial Dehn filling along $K$.\endproclaim 

{\bf Proof:} For basic concepts on Dehn 
filling and Dehn surgery, the reader is referred to [Ro].

\ssk

To show that \clsp is essential in the manifold $K(r)$ obtained by 
$r$-Dehn-filling along $K$, for $r\neq$1/0, we will verify the six properties
(1)-(6) above, in this new setting.
Once again, the first four of these properties require no extra proof, 
since $N(B)$ has not changed; only where it is embedded has. So we need only show that 
$K(r)$\smin int$N(B)$ = \emrsp is irreducible, and \delh $N(B)$ is 
incompressible in \emr. 

We can think of \emrsp as the result of $r$-Dehn-filling on $K$ in 
the genus-2 handlebody \eno. Both of our proofs will rely on the fact that 
\eno\smin int$(N(K))$ contains two embedded annuli $A_i$ =$D_i$\smin int$(N(K))$, 
$i$=1,2, each with one component on \delh $N(B)$ and the other a meridional 
loop on \del $N(K)$. 

To show that \emrsp is irreducible, suppose it is not. Then there is a 
2-sphere $S$ in \emrsp which does not bound a ball in \emr. Choose such a 
sphere which intersects (transversely) the (image of the) knot $K$ in the 
fewest number of points. It is then standard that $S$\smin int$(N(K))$ = $S^\prime$ 
is an incompressible and \del-incompressible 
planar surface in \eno\smin int$(N(K))$ = \emo. 
The curves $S^\prime$/ap\del$N(K)$ are parallel curves of slope $r$.

Look at $S^\prime$\ap$A_i$; it consists of circles and arcs. Trivial circles 
of intersection can be removed by isotopy, since $S^\prime$ 
is incompressible. The arcs of intersection cannot meet 
the boundary component of $A_i$ coming from \del\eno, since $S$ misses \del\emr. 
These arcs of intersection are therefore boundary parallel, and so can also be 
removed by isotopy, since $S^\prime$ is \del-incompressible. After these 
isotopies, if \del$S^\prime\neq\emptyset$, \del$S^\prime$\sset\del$N(K)$ misses 
a meridional loop, and hence consists of meridional loops. So $r$ = 1/0, a 
contradiction.

Therefore \del$S^\prime$=$\emptyset$, i.i., $S^\prime$ = $S$. But since 
\enosp is irreducible, $S$ bounds a 3-ball in \eno. This 3-ball must intersect, 
hence contain, $K$, since otherwise $S$ bounds a 3-ball in \emr. But this 
implies that $K$ is null-homologous in \eno, which is impossible since it 
intersects a compressing disk $D_1$ of \del\eno\ exactly once. So the reducing sphere 
cannot exist; \emrsp is irreducible.

\ssk

To prove incompressibility of \delh, we again appeal to the two disjoint 
annuli $A_i$ in \emo, joining meridional loops in \del $N(K)$ to loops in 
\delh $N(B)$. The two loops  in \delh $N(B)$ are obviously not homotopic to 
one another on \delh $N(B)$; they lie in different components. We have already 
seen above that \delh $N(B)$ is incompressible in \eno\smin int$(N(K))$. 
It then follows from Theorem 4 of [Me] that \delh $N(B)$ will 
remain incompressible in any manifold obtained by non-trivial (i.e., 
non-meridional) Dehn filling along $K$ (in \eno). 

\ssk

Therefore, all of the properties of essentiality for $B$ (and hence for \cl) are 
satisfied in any manifold obtained by non-trivial Dehn filling along any of 
the knots $K$. So $K$ is persistently laminar. \,\, \qed

\ssk

\heading{\S 4 \\ Generalizations}\endheading

The existence of a persistent lamination in the complement of knots $K$ 
obtained from $T_0$ duplicates previous work. If $T$ is a rational tangle, 
then $K$ is a Montesinos knot, and for such tangles Delman [De] has constructed 
persistent laminations for the resulting knots. On the other hand, 
if $T$ is an non-split tangle, then Wu [Wu] has shown that $K$ 
admits a persistent lamination. These two results have very powerful 
generalizations, as well. The intersection of the complements of these 
two classes of tangles is the collection of split, non-rational, tangles, 
and so the resulting knots are all connected sums with a square knot. One 
of the swallow-follow tori for each knot will then remain incompressible 
under all non-trivial Dehn fillings. 

The technique for building the branched surface $B$ that we have used here 
can, however, be easily extended to more than one tube; see Figure 7 for 
the case $n$=2. We can then string arcs through the tubes of $B_n$, to 
create a tangle of 2n arcs in a 3-ball. $S^3$\smin int$(N(B_n))$ = \ennsp 
is again a handlebody (it is a 3-ball with 2$n$ 1-handles attached,) and 
the branch curve \gamsp separates the 1-handles into two collections of 
$n$ each (Figure 8). In this case however, we cannot always add arcs 
in the central 3-ball to create a knot in any way we choose, and still 
expect $B_n$ to be essential; an arc running (parallel to \del\enn) from 
arcs on the same side of \gamsp, for example, has a (parallel) \del-compressing 
disk around it (see Figure 8b).

\ssk
\leavevmode

\epsfxsize=5in
\centerline{{\epsfbox{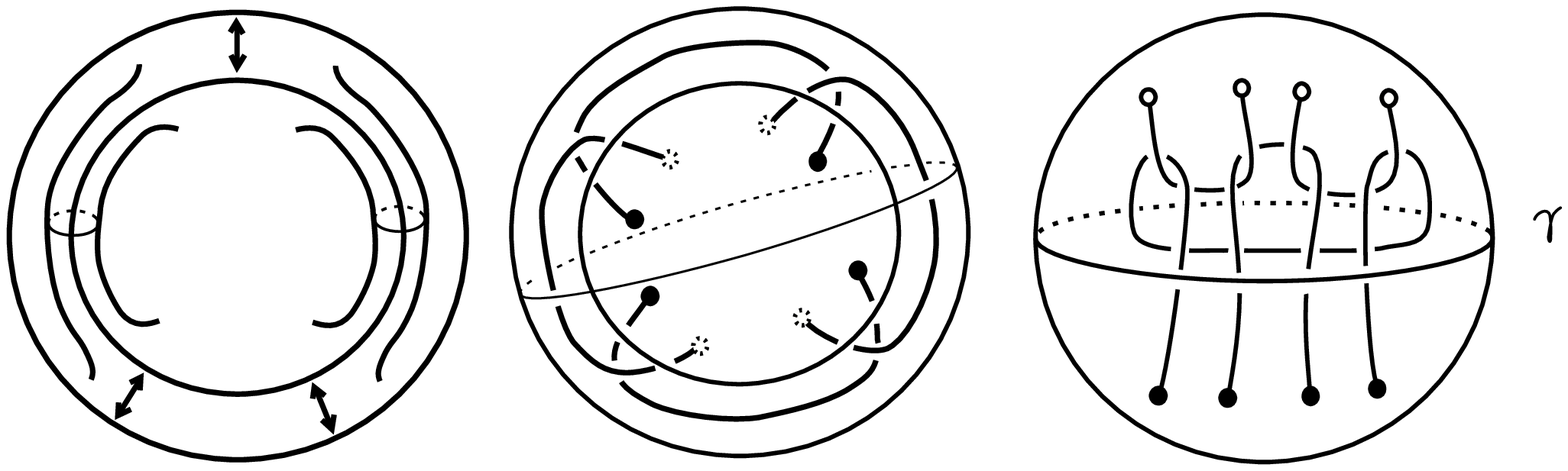}}}

\centerline{Figure 7}

\ssk
\leavevmode

\epsfxsize=5in
\centerline{{\epsfbox{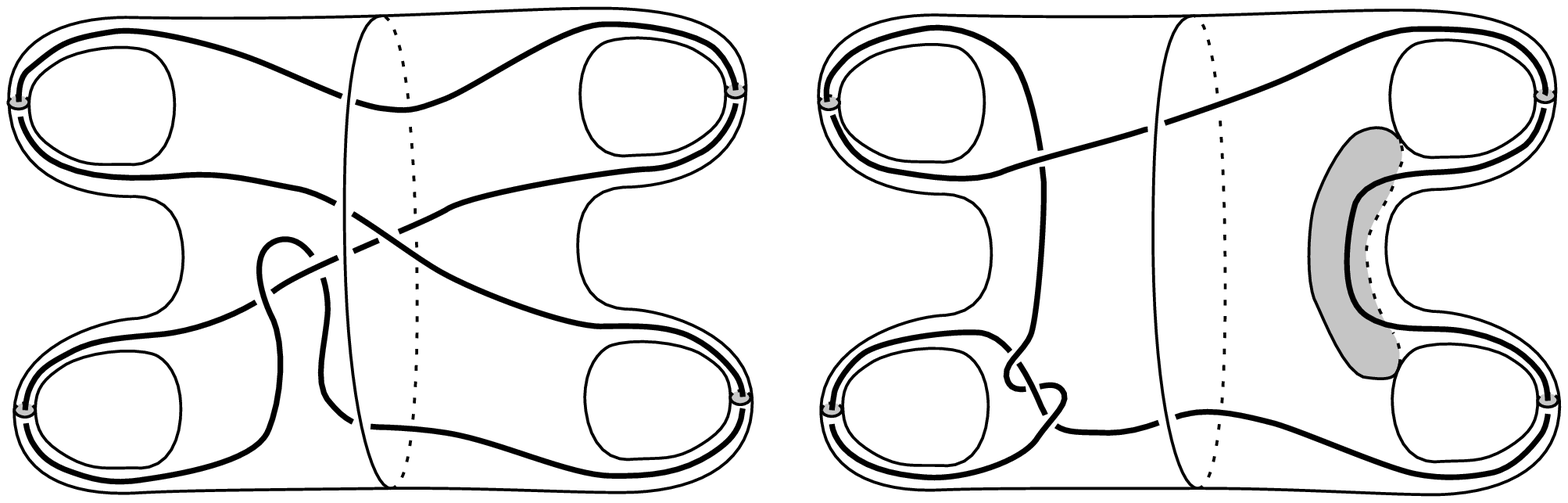}}}

\centerline{Figure 8}

\ssk

However, incompressibility of the horizontal boundary \delh$N(B_n)$ in \ennsp 
is the only obstruction to the essentiality of $B_n$, as well as remaining 
essential under any non-trivial Dehn surgery. The first four conditions on 
essentiality follow the exact same line as in our original case. We also still 
have the meridional annuli which allow us to verify that irreducibility and 
incompressibility of \delh$N(B_n)$ will be inherited under any non-trivial Dehn 
filling. In some sense, it turns out, the phenomenon described above is also 
the only way to prevent incompressibility, as well.

We can push compressing disks off of the tubes, by pushing them off of the 
meridional annuli; see Figure 9. Trivial circles of intersection with $A_i$ 
can be removed by isotopy, since \ennsp (and therefore \enn\smin $K$) is 
irreducible, and then we may surger along arcs of intersection to create two 
disks, at least one of which has boundary non-trivial in \delh $N(B)$, giving 
us a new compressing disk with fewer intersections with $A_i$. Finally, we 
cannot have any circles of intersection which are essential in $A_i$, since 
surgering along the innermost one (using the disk in $D_i$ that it bounds, which meets $K$ 
once) would produce a disk and a 2-sphere (in $S^3$) each intersecting $K$ exactly once.
But a sphere in $S^3$ cannot meet a knot only once.

\ssk
\leavevmode

\epsfxsize=5in
\centerline{{\epsfbox{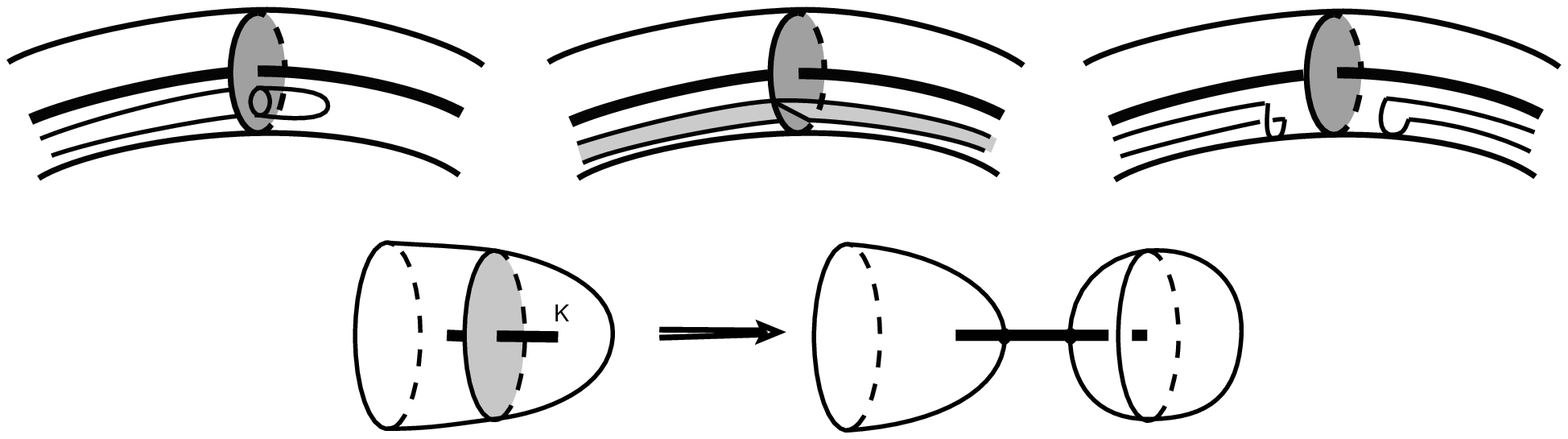}}}

\centerline{Figure 9}

\ssk

Our compressing disk $D$ then lies in the central 3-ball 
piece $B^3$ of $S^3$\smin int$(N(B_n))$. Again, we can assume (by disk-swapping) 
that $D$ misses the obvious compressing disk $D^\prime$ that the branch curve 
\gamsp bounds. $D$ then splits $B^3$ into two 3-balls $B^3_1$ and $B^3_2$; one 
of them, $B^3_1$ say, misses $D^\prime$. We must then have 
$K$\ap$B^3_1\neq\emptyset$; otherwise, $D$ can be isotoped, rel boundary, into 
\delh$N(B_n)$, since \del$B^3_1$ misses $K$, hence misses the subdisks of 
\del$B^3$ which the 1-handles of \ennsp are attached to (see Figure 10).

Therefore $K$\ap$B^3_1$ consists of some non-zero number of components of the 
2$n$-strand tangle $K$\ap$B^3$. These arcs are disjoint from $B^3_2$, and so 
are disjoint from $D^\prime$, and so each joins endpoints of core arcs of 
1-handles which are on the \ubr{same} side (i.e., the $D$-side) of \gam.

\ssk

\leavevmode

\epsfxsize=1.5in
\centerline{{\epsfbox{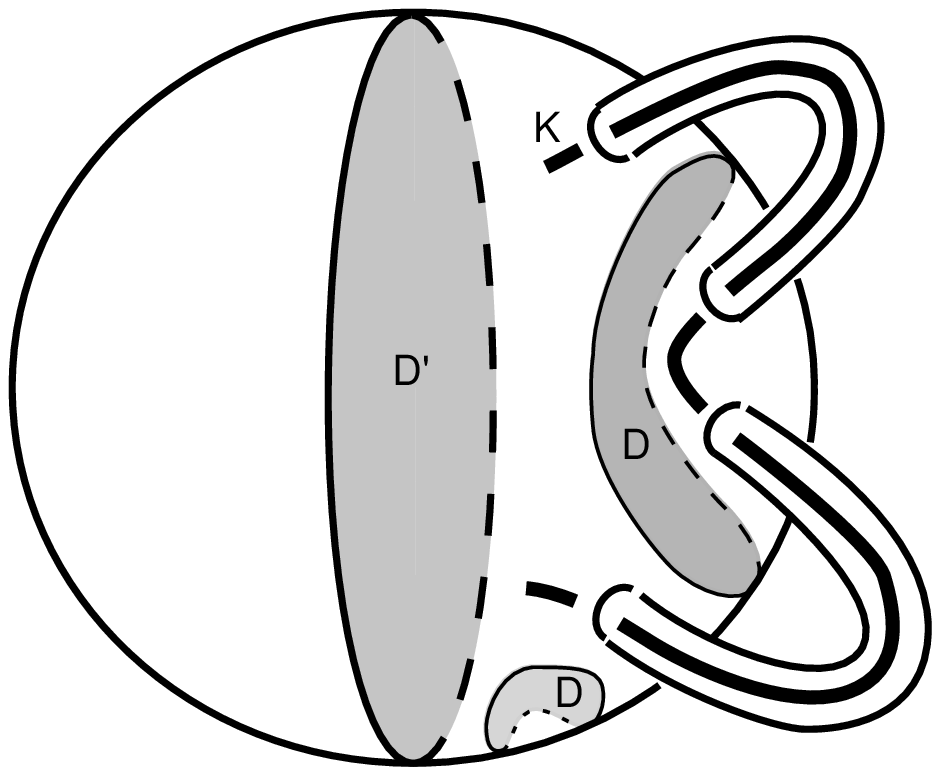}}}

\centerline{Figure 10}

\ssk

\ubr{Therefore}, one way to ensure that the branched surface $B_n$ is essential 
in the complement of $K$, and remains esential under all non-trivial Dehn 
fillings, is to insist that all of the arcs we use to build $K$ travel from 
one side of \gamsp to the other, as in Figure 8a. This is still a vast number 
of knots, all of which admit a persistent lamination. Even more, we can allow 
ourselves to connect the ends of this 2$n$-strand tangle as above to create 
links, as well. Since each component of the link must visit both sides of \gam, 
each component comes equipped with two of the meridional annuli $A_i$, whose 
boundary components are on different components of \delh $N(B)$. Therefore, 
non-trivial Dehn filling on each component of the link (what is sometimes 
called a {\it complete} Dehn filling) yields a manifold in which our lamination 
remains essential. Note that, with our original tangle $T_0$, the condition 
that the arcs of the tangle $T$ travel from one side of \gamsp to the other is 
precisely the condition that the resulting link is in fact a knot. So this new 
condition is a natural extension.

\ssk

\leavevmode

\epsfxsize=3in
\centerline{{\epsfbox{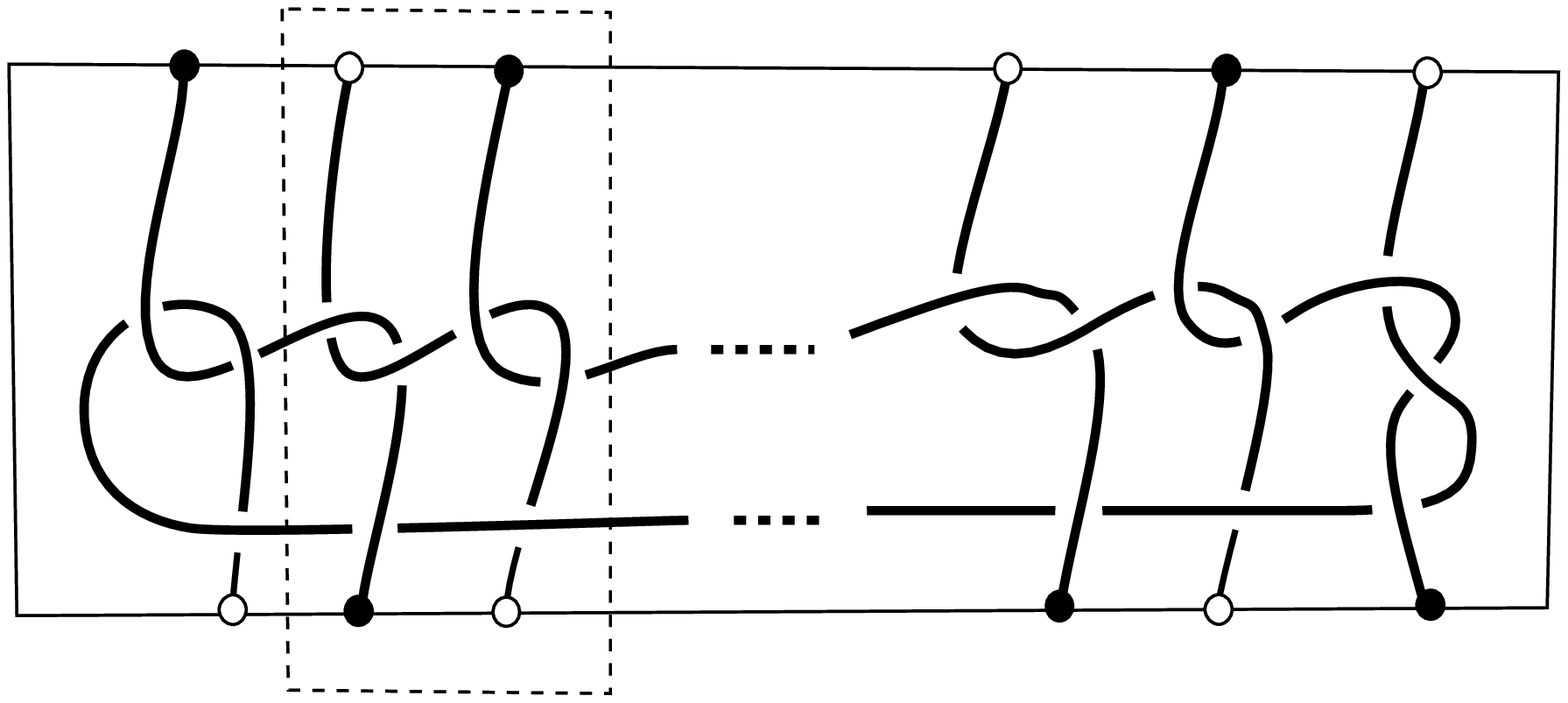}}}

\centerline{Figure 11}

\ssk

It is not hard to see that the tangle of Figure 7 (and its generalizations
with more tubes) 
can be isotoped to an alternating diagram. Adding more arcs to the tangle 
(i.e., adding more tubes to $B_n$) simply amounts to grafting on an additional 
fundamental piece to the tangle, shown in the dotted rectangle in Figure 11. 
The condition above then amounts to requiring that the complementary 2$n$-strand 
tangle join black-dotted ends to white-dotted ends. Joining together all but two 
pairs of ends produces an ordinary tangle. Note that each strand of the tangle 
must have been built from an odd number of our original strands, in order for 
its ends to lie on the same side of the compressing disk $D$. Since our 
lamination remains essential and persistent no matter how this tangle is 
completed to a knot, these tangles are persistently laminar.

\vfill\eject

\heading{\S 5 \\ Still more generalizations}\endheading

\ssk

\leavevmode

\epsfxsize=4in
\centerline{{\epsfbox{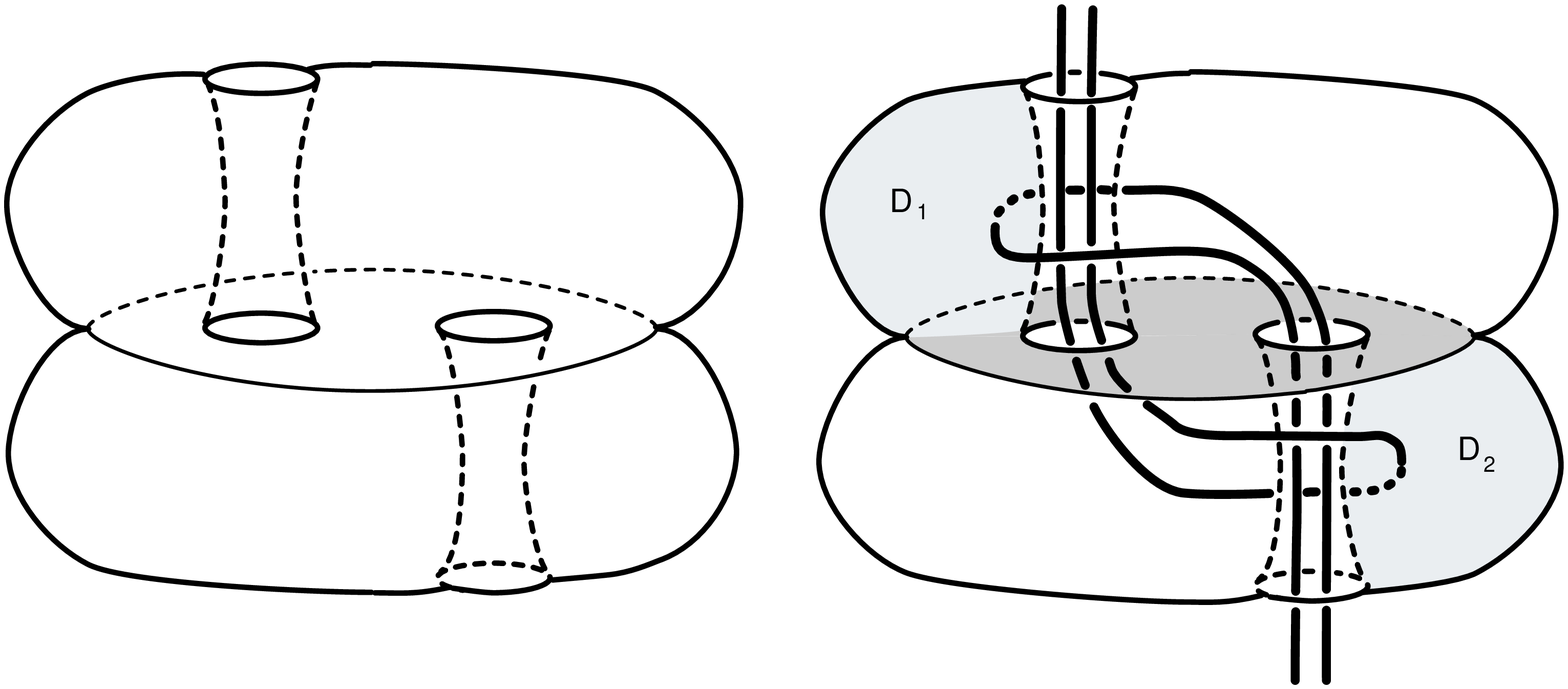}}}

\centerline{Figure 12a \hskip2in Figure 12b}

\ssk

Ramin Naimi has pointed out that the branched surface $B$ which we began with 
can be drawn in a different (and ultimately more useful) way; see Figure 12. 
In this form it is easy to see all of the components of the construction 
which we have exploited; the compressing disk bounded by the branch curve 
\gam, the compressing disks for the two 1-handles, and the tangle $T_0$ 
built from the core arcs of the 1-handles.

\ssk

\leavevmode

\epsfxsize=4in
\centerline{{\epsfbox{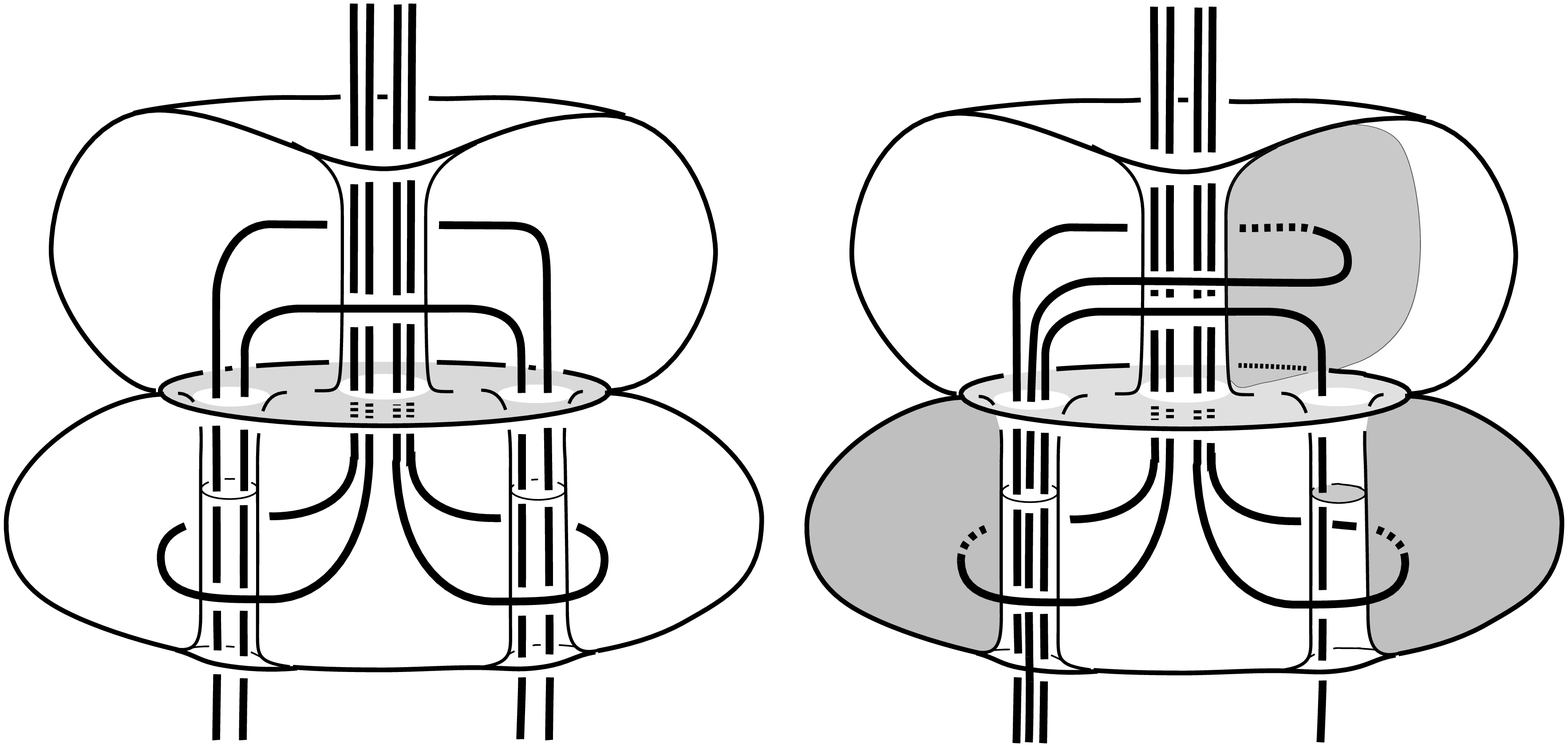}}}

\centerline{Figure 13}

\ssk

Our more general branched surfaces $B_n$ have similarly simple pictures; 
see Figure 13a. In this form, however, it is also easy to see that there 
are \ubr{different} choices of how to write $S^3$\smin int$B_n$=\ennsp 
as a 3-ball $B^3$ (containing the branch curve \gam) with 1-handles 
attached, by choosing \ubr{different} compressing disks for \del\enn\smin\gam; 
see Figure 13b. We can take the core arcs of these compressing disks, and 
think of them as a 2$n$-strand tangle in the complementary 3-ball 
$S^3$\smin$B^3$. Structurally, these tangles have the exact same properties 
which we used in Section 4 to show that $B_n$ is essential in the 
complement of any knot or link obtained by gluing on a tangle in $B^3$ all 
of whose strands cross the disk $D$ bounded by \gam. Therefore, we can 
obtain new examples of persistently laminar tangles by choosing sets of 
compressing disks for the two genus-$n$ handlebodies of \enn$|D$, and taking 
the core arcs of the disks. We give a further example in Figure 14.

\ssk

\leavevmode

\epsfxsize=2in
\centerline{{\epsfbox{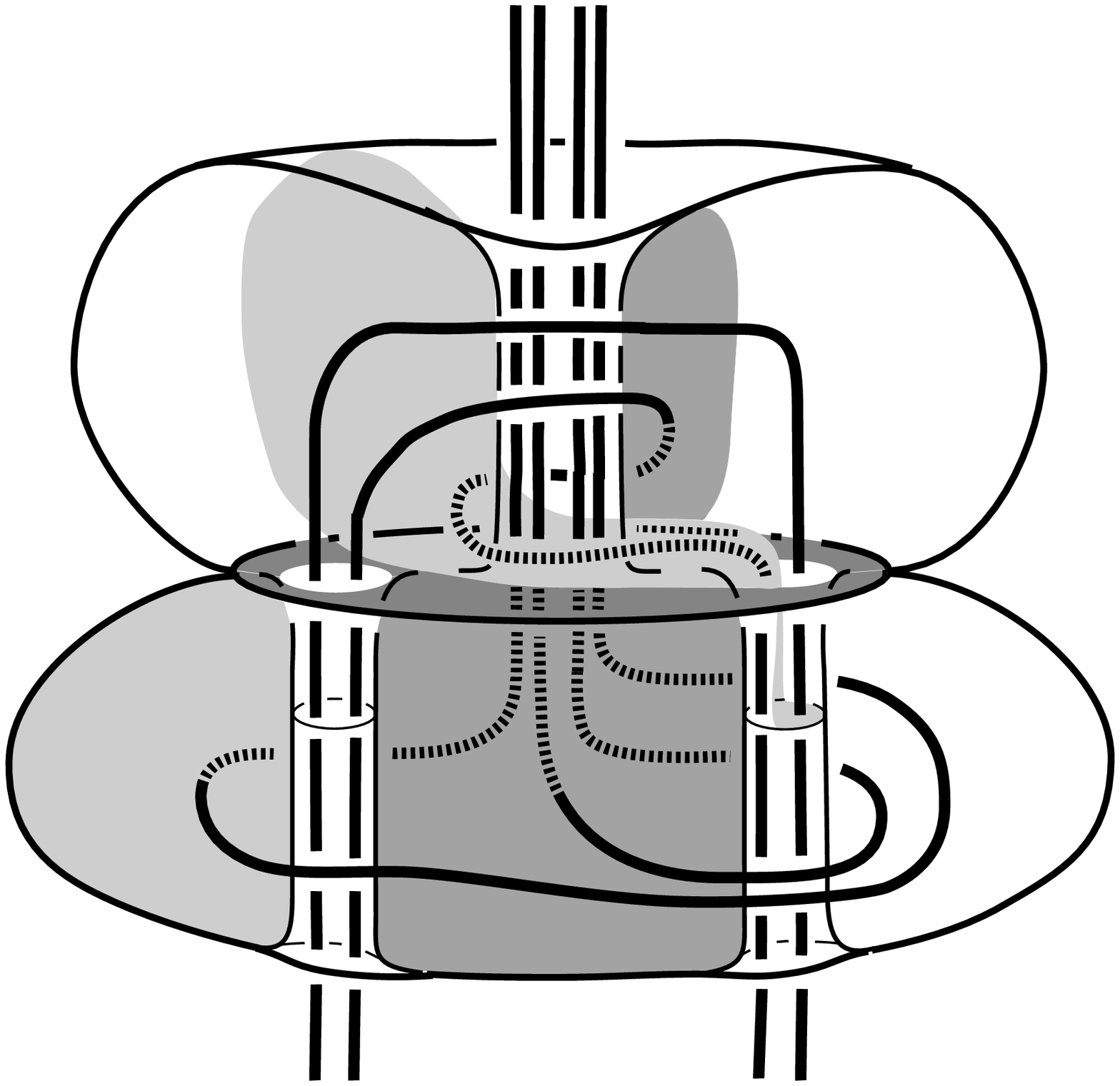}}}

\centerline{Figure 14}

\ssk

We can list of the properties of our branched surfaces $B_n$ 
which we have used in our proofs; this gives us a recipe for finding 
persistently laminar tangles. 
We needed a transversely orientable branched surface $B$ in $S^3$ having 
one branched curve \gam, with no triple points, so that $B$\smin\gamsp is 
connected. We also require that \gamsp bounds a disk $D$ which splits 
$S^3$\smin int$N(B)$ into two genus-$n$ handlebodies (note that the two 
handlebodies must have the same genera); in particular, $S^3$\smin int$N(B)$ 
is a handlebody. Choosing compressing disks for each handlebody and taking 
their core arcs gives us a 2$n$-strand tangle which we can string together 
as above to create persistently laminar tangles. In Figure 15 we
provide an example, using this recipe.

\ssk
\leavevmode

\epsfxsize=3in
\centerline{{\epsfbox{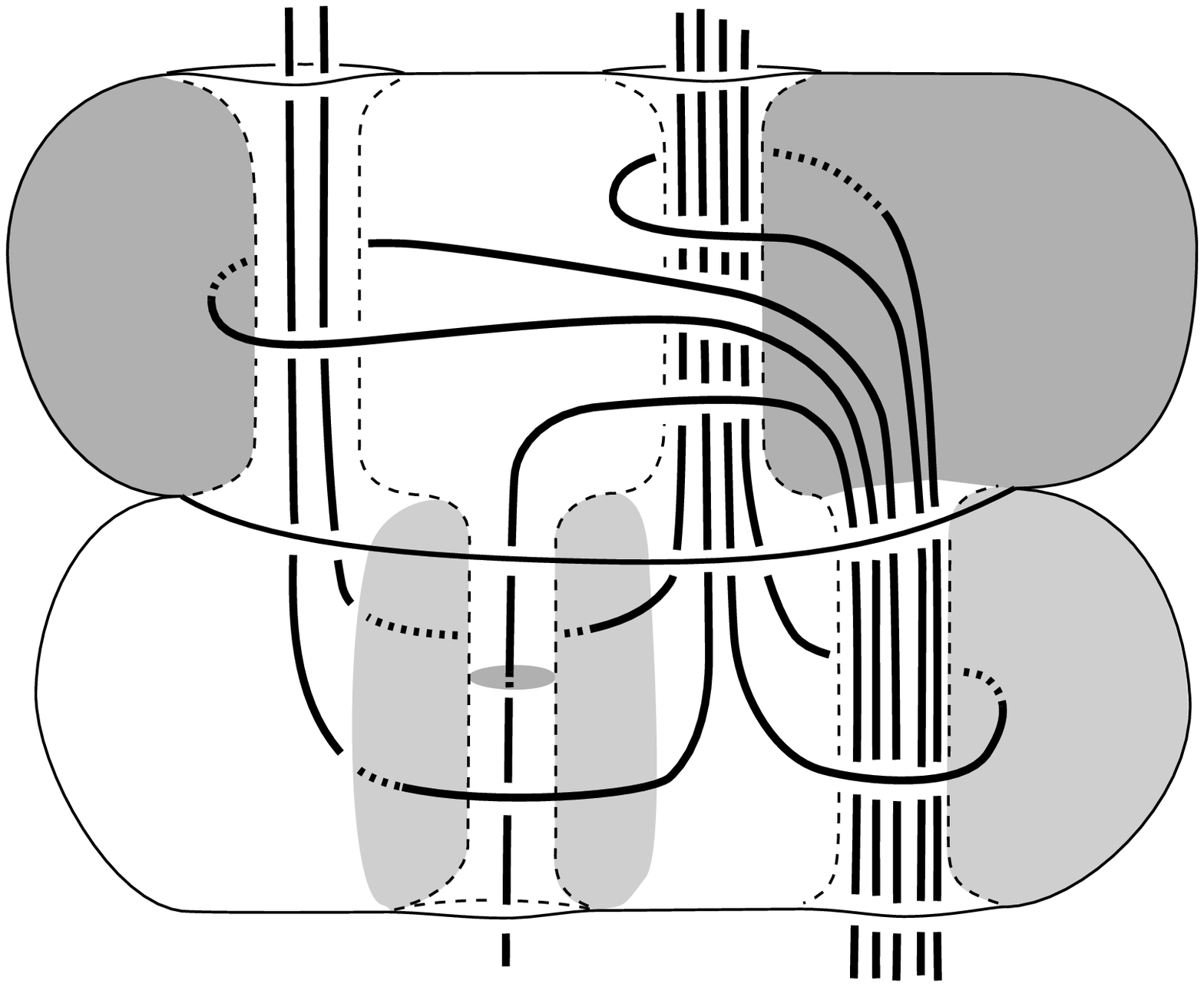}}}

\centerline{Figure 15}

\heading{\S 6 \\ Concluding remarks}\endheading

The tangles we have described here come equipped with a lamination in their 
tangle space (the 3-ball with the arcs of the tangle removed), which remains 
essential after non-trivial surgery on any knot constructed from the tangle. 
This is what we have called persistently laminar. One could weaken this 
definition, without losing its essential strength, by requiring instead 
that for every knot obtained from the tangle, there is a lamination which 
remains essential under non-trival Dehn filling. The work of Delman [De] 
and Wu [Wu] then demonstrate that many other tangles are persistently 
laminar in this sense; for example, the sum of two rational tangles 
whose associated rational numbers have denomentators at least 3 and have different 
signs (such as, for example, our tangle $T_0$) [De], or the sum of two 
atoroidal tangles [Wu]. Most algebraic tangles (see [Wu]) are also persistently 
laminar, in this weaker sense. The technique of the previous paragraph can 
easily provide examples of persistently laminar (alternating) tangles which 
cannot be decomposed (non-trivially) as the sum of two tangles, however, 
making them disjoint from these collections of tangles.

There are, of course, many tangles which are \ubr{not} persistently laminar; 
any tangle which can be summed to give a knot admitting a finite or reducible 
surgery, for example, cannot be persistently laminar, because the
surgery manifold is not laminar. So, for example, no 
rational tangle is persistently laminar; each can be summed with another 
rational tangle to produce a (2,q)-torus knot. Other, more sporadic, examples 
can easily be given.

\ssk

In this paper we have worked, in some sense, backwards, by building a 
lamination and then finding the tangle space which it should live in. A far 
more difficult (and so correspondingly rewarding) approach is to try to 
determine if a given tangle is persistently laminar, in either sense. Wu, 
for example, suggests the tangle of Figure 16 as an example; it is, in some sense, the 
smallest non-algebraic tangle. We do not know whether or not it is persistently laminar. 
No knot obtained from it by tangle sum with another tangle is known to fail to
be persistently laminar.

\ssk
\leavevmode

\epsfxsize=2in
\centerline{{\epsfbox{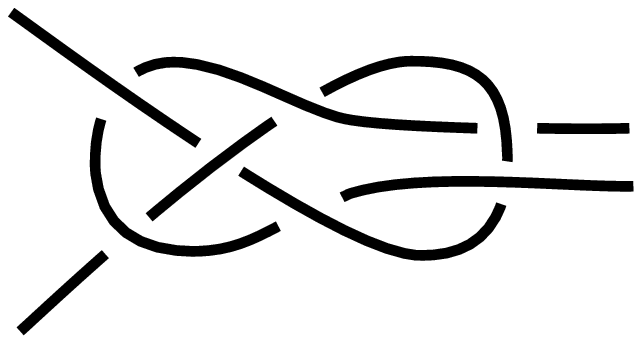}}}

\centerline{Figure 16}

\ssk

The laminations we have worked with are also in some sense the 
`simplest' laminations one could build; their branched surfacees have a 
single branch curve with no triple points. It is remarkable how many knots 
these very simple laminations are persistent for; are there other constructions 
which are similarly powerful?

\smallskip

\centerline{\bf Acknowledgements}

\smallskip

The author is very grateful to Jim Hoffman, Ramin 
Naimi, Ulrich Oertel, Rachel Roberts, and Ying-Qing Wu for many helpful 
conversations during the preparation of this work. This research was 
supported in part by NSF grants DMS$-$9400651 and DMS$-$9704811.

\end